\documentclass{IEEEtaes}

\jvol{XX}
\jnum{XX}
\jmonth{XXXXX}
\paper{1234567}
\pubyear{2024}
\doiinfo{TAES.2024.Doi Number}

\begin{document}

\title{Recursive Polynomial Method for Fast Collision Avoidance Maneuver Design} 

\author{ZENO PAVANELLO, LAURA PIROVANO and \\ ROBERTO ARMELLIN}
\affil{Te Punaha Atea - Space Institute, The University of Auckland, Auckland, New Zealand} 

% \author{LAURA PIROVANO}
% \affil{Te Punaha Atea - Space Institute, The University of Auckland, Auckland, New Zealand} 

% \author{ROBERTO ARMELLIN}
% % \member{Member, IEEE}
% \affil{Te Punaha Atea - Space Institute, The University of Auckland, Auckland, New Zealand}

\receiveddate{Manuscript received May 16, 2024; revised XXXXX 00, 0000; accepted XXXXX 00, 0000.\\
© 2024 IEEE. Personal use of this material is permitted. Permission from IEEE must be obtained for all other uses, in any current or future media, including reprinting/republishing this material for advertising or promotional purposes, creating new collective works, for resale or redistribution to servers or lists, or reuse of any copyrighted component of this work in other works. \\
This work has been submitted to the IEEE for possible publication. Copyright may be transferred without notice, after which this version may no longer be accessible. \\
The codes related to this work can be downloaded for reproducibility at \url{https://github.com/zenop95/Recursive-polynomial-CAM-design/releases/tag/v1.1.0}.\\
}
%% \accepteddate{XXXXX XX XXXX}
%% \publisheddate{XXXXX XX XXXX}

\corresp{\itshape{(Corresponding author: Z. Pavanello).}}

% \authoraddress{Zeno Pavanello is a Ph.D. Student at Te P\=unaha \=Atea - Space Institute, The University of Auckland, Auckland, 1010 New Zealand
% (e-mail: \href{mailto:zpav176@aucklanduni.ac.nz}{zpav176@aucklanduni.ac.nz}.\\
% Laura Pirovano is a Postdoctoral Researcher at Te P\=unaha \=Atea - Space Institute, The University of Auckland, Auckland, 1010 New Zealand
% (e-mail: \href{mailto:laura.pirovano@auckland.ac.nz}{laura.pirovano@auckland.ac.nz}).\\
% Roberto Armellin is a Professor at Te P\=unaha \=Atea - Space Institute, The University of Auckland, Auckland, 1010 New Zealand
% (email: \href{mailto:roberto.armellin@auckland.ac.nz}{roberto.armellin@auckland.ac.nz}).}

\markboth{Z. PAVANELLO ET AL.}{RECURSIVE POLYNOMIAL METHOD FOR FAST COLLISION AVOIDANCE MANEUVER DESIGN}
\maketitle

\begin{abstract}
A simple and reliable algorithm for collision avoidance maneuvers (CAMs), capable of computing impulsive, multi-impulsive, and low-thrust maneuvers, is proposed. The \acrfull{poc} is approximated by a polynomial of arbitrary order as a function of the control, transforming the \acrshort{cam} design into a \acrlong{pp}. The solution procedure is initiated by computing the \acrshort{cam} via a first-order greedy optimization approach, wherein the control action is applied in the direction of the gradient of \acrshort{poc} to maximize its change. Successively, the polynomial is truncated at higher orders, and the solution of the previous order is used to linearize the constraint. This enables achieving accurate solutions even for highly nonlinear safety metrics and dynamics. 
Since the optimization process comprises only polynomial evaluations, the method is computationally efficient, with run times typically below $1$ \si{s}. Moreover, no restrictions on the considered dynamics are necessary; therefore, results are shown for Keplerian, J2, and \acrlong{cr3bp} dynamics.
\end{abstract}


\section{INTRODUCTION}
According to the most recent European Space Agency (ESA) Space Debris report \cite{ESA2023}, the number of cataloged \acrlong{rso}s is over 32,000. This exponential increase, fueled by spacecraft miniaturization and the deployment of mega-constellations like Starlink, is reshaping how we think about space traffic management. A cluttered environment amplifies the frequency of conjunctions and drives the necessity for a commensurate increase in the number of \glspl{cam}.

Presently, conjunction analysis and collision avoidance operations are predominantly executed by operators on the ground, leveraging tools and processes refined over the past two decades. While these tools facilitate operational activities, the reliance on human intervention in decision-making processes and maneuvers design will become unsustainable. Therefore, the demand for automated conjunction screening, \gls{ca} decision-making, and \gls{cam} design and execution intensifies. This calls for \gls{cam} algorithms suitable for autonomous and, potentially, onboard computations providing fuel-optimal maneuvers in a short execution time. This paper aims to address this pressing need by proposing a method for \gls{cam} optimization suitable for autonomous applications. 

A \gls{cam} is performed when operators receive a \gls{cdm} indicating the \gls{tca} between the satellite and a secondary object and the relevant information of the conjunction. Generally, the estimated \gls{poc} at \gls{tca} exceeds a threshold value, and the objective of the \gls{cam} is to lower it while minimizing the total $\dv$. Alternatively, if considerations on the uncertainty of the state of the two objects are not of interest, the conjunction is mitigated by increasing the miss distance at \gls{tca}.
Previous research has explored various methodologies for optimizing \glspl{cam}, often simplifying the problem by assuming small maneuvers \cite{Patera2003}.
Alfano \cite{Alfano2005Collision} introduced a \gls{cam} analysis tool capable of conducting parametric studies on single-axis and dual-axis maneuvers, while the German Aerospace Center and ESA have developed their own CAM optimization tools with differing degrees of flexibility and optimality \cite{Aida2016,GrandeOlalla2013}. Moreover, recent advancements have seen the emergence of multi-objective approaches for CAM design, enabling comprehensive analyses and the identification of Pareto optimal solutions \cite{Morselli2014Collision}. Even though very useful for robust and optimal maneuvers computations, all of these solutions are not feasible for autonomous applications because of their high computational load. In recent years, a substantial corpus of research has emerged to devise efficient \gls{cam} algorithms. Bombardelli \cite{Bombardelli2014} and Bombardelli and Hernando-Ayuso \cite{Hernando-Ayuso2020,Bombardelli2015}  devised analytical and semi-analytical methods to minimize miss distance or collision probability using a single impulse for a given magnitude of $\dv$, demonstrating promising convergence properties. De Vittori \textit{et al.}  proposed an analytical approach for the energy-optimal problem and a semi-analytical one for the fuel-optimal to handle low-thrust \glspl{cam} \cite{DeVittori2022}, which are particularly relevant given the increasing popularity of electric thrusters \cite{Mazouffre2016}. Direct optimization methods have also been extensively utilized to compute \glspl{cam}. Misra and Dutta \cite{Dutta2022} linearize the \gls{poc} constraint and achieve fast but conservative solutions, so they do not guarantee the minimum fuel expenditure. Armellin \cite{Armellin2021} convexifies the fuel-optimal \gls{cam} optimization \gls{nlp} using \gls{da}, lossless convexification, and transforming the \gls{poc} constraint into a keep-out-zone constraint; the convex program is then solved iteratively, and the original fuel-optimal solution is achieved. Armellin's method was later extended to solving long-term encounters \cite{Pavanello2023Long}, multiple consecutive encounters \cite{Pavanello2024Multiple}, and encounters during low-thrust arcs \cite{Pavanello2024LowThrust}.

In this work, we propose an alternative method to design the \gls{cam}, which can retain computational speed comparable to semianalytical indirect methods while encompassing flexibility similar to direct methods. We frame the \gls{cam} design problem as a \gls{pp} by utilizing the energy-optimal objective function and approximating \gls{poc} as a high-order Taylor polynomial in the control. This can be represented as a series of impulses or accelerations applied at specified nodes or segments. The \gls{pp} is solved by leveraging the assumption that high-order terms get smaller in a convergent Taylor series.
Thus, starting from a greedy solution of the problem with a linearized \gls{poc} constraint, we establish a simple and computationally efficient recursive approach that can achieve quasi-optimal results while accurately meeting the safety constraint. The \gls{pp} is solved using a succession of solutions of increasing polynomial order. The tensor notation, notably used in astrodynamics problems in recent years to define state transition tensors \cite{Park2006}, is adopted to frame the recursive method. \gls{da} is used to compute the tensors efficiently. This methodology allows for considering any dynamics, including high-accuracy Earth orbit models and \gls{cr3bp} models for the Cislunar environment. Given its efficiency, flexibility, and simplicity, this tool could become an important asset for operators to transition towards autonomous \glspl{cam} design.

The solution obtained with the proposed recursive method is compared with the interior-point method implemented in MATLAB's fmincon for the fuel-optimal \gls{nlp}. In this way, the recursive method is proven to be faster than a state-of-the-art solver with limited loss in terms of optimality.

The paper is organized as follows. In \cref{sec:dyn}, the dynamics models considered for the \gls{cam} applications are introduced. \cref{sec:optProb} sets the \gls{nlp} using automatic Taylor expansions. The main innovation proposed by this work, i.e., the recursive algorithm to solve the \gls{cam} \gls{pp}, is shown in \cref{sec:recursive}. In \cref{sec:results}, operational results are presented, and the comparison with the interior-point solver is discussed. Lastly, in \cref{sec:conclusions}, conclusions are drawn. The open-source software developed for this work can be consulted and downloaded using the link on the first page.

\section{DYNAMICS FRAMEWORKS}
\label{sec:dyn}
This section introduces the framework of the \gls{cam} optimization problem dynamics. Different models are considered, depending on the domain in which the \gls{cam} needs to be performed and on the degree of accuracy of the modeled environment.

\subsection{Conjunction in the B-plane}

This work assumes that the relative velocity between the two satellites at \gls{tca} is high. Therefore, the encounter can be regarded as an instantaneous event \cite{Armellin2021}. 
A \gls{cdm} gives the state and covariance matrix of the two satellites involved in the conjunction at \gls{tca}. By international standards, a \gls{cdm} assumes that the states of the two bodies are Gaussian multivariate random variables, so the mean state and covariance are sufficient to describe them. Typically, for Earth orbit conjunctions, the states are expressed in an \gls{eci} reference frame; in the Cislunar environment, they are expressed in the synodic reference frame. Therefore, the mean state of the primary is $\vec{x}(t_{CA})=[\vec{r}(t_{CA}); \hspace{3pt} \vec{v}(t_{CA})]$, while the debris has mean state $\vec{x}_s(t_{CA})=[\vec{r}_s(t_{CA}); \hspace{3pt} \vec{v}_s(t_{CA})]$. Here, $\vec{r}(t_{CA}),\vec{r}_s(t_{CA}),\vec{v}(t_{CA})$ and $\vec{v}_s(t_{CA})\in \mathbb{R}^3$ denote the position and velocity of the centers of mass of the objects. Analogously, the covariances of the two objects are $\vec{C}$ and $\vec{C}_s\in \mathbb{R}^{6\times6}$. The relative state is defined as the subtraction of the two absolute states. Since the relative state is of interest only at \gls{tca}, we omit the argument
\begin{equation}
\begin{aligned}
    \vec{x}_{rel} = \vec{x}(t_{CA}) - \vec{x}_s(t_{CA}), \quad & \vec{C}_{rel} = \vec{C} + \Vec{C}_s.
\end{aligned}
\end{equation}
Note that $\vec{x}_{rel}=[\vec{r}_{rel}; \hspace{3pt} \vec{v}_{rel}]$ and $\vec{C}_{rel}=[\vec{C}_{rr} \hspace{3pt} \vec{C}_{rv}; \hspace{3pt} \vec{C}_{rv} \hspace{3pt} \vec{C}_{vv}]$. For the scope of this study, we are only interested in the positional part of the covariance matrix, which will be referred to as $\vec{P} = \vec{C}_{rr}$ in the following.

Given the short-term nature of the encounter, it can be studied in the B-plane reference frame denoted as $\mathcal{B}$. %As depicted in \cref{fig:B-plane},
$\mathcal{B}$ is centered on the secondary object;
the $\eta$ axis aligns with the direction of the relative velocity of the primary with respect to the secondary, while the $\xi\zeta$ plane is perpendicular to the relative velocity axis. Since the \gls{tca} is the moment when the miss distance is minimum, $\vec{r}_{rel}\cdot\vec{v}_{rel} = 0$ by assumption, indicating that the relative position lies within the $\xi\zeta$ plane. 
Therefore, for the computation of \gls{poc}, we will employ the projection of the relative state and its covariance on the $\xi\zeta$ plane: $\Vec{r}_\mathcal{B}\in\mathbb{R}^2$ and $\Vec{P}_\mathcal{B}\in\mathbb{R}^{2\times2}$.
The subsequent discussion assumes that all uncertainty is concentrated around the secondary object, while all mass is concentrated around the primary \cite{Li2022}.

% \begin{figure}[tb!]
%     \centering
%     \input{Figures/B-plane}
%     \vspace{.5cm} 
%     \caption{B-plane construction (adapted from \cite{Armellin2021}).}
%     \label{fig:B-plane}
% \end{figure}

A conservative way to compute \gls{poc} of the encounter is enveloping each of the bodies in a sphere with a radius equal to the largest dimension. The sum of the radii of the two spheres is the \gls{hbr}, which defines the combined hard body sphere centered in the primary. The projection of this sphere onto the B-plane is the combined hard body circle $\mathbb{C}_\mathrm{HBR}$.
To compute \gls{poc}, denoted as $P_C\in\mathbb{R}_\mathrm{+}$, the probability density function of the projection of the relative position onto the B-plane is integrated over the combined hard-body circle, as follows
\begin{equation}
\begin{split}
\pc = & \frac{1}{(2 \pi)^{3/2}\sqrt{\mathrm{det}\big(\boldsymbol{P}_\mathcal{B}\big)}} \cdot \\& \cdot \iint_{\mathbb{C}_\mathrm{HBR}} \mathrm{exp}\left(-\frac{\vec{r}_\mathcal{B}\transp\vec{P}_\mathcal{B}^{-1}\vec{r}_\mathcal{B}}{2}\right) \mathrm{d} A.
\end{split}
\label{eq:pceq}
\end{equation} 

Multiple approaches have been proposed to approximate the integral in \cref{eq:pceq} \cite{Li2022}. The results presented in \cref{sec:results} are obtained using Chan's method \cite{Chan2004International}, which approximates \gls{poc} through a convergent series that involves equivalent cross-sections. Nonetheless, other methods are present in the literature, and we highlight the fact that our optimization method is agnostic to the PoC model used since the automatic Taylor expansion can deal with any function. Therefore, the following discussion is made without assuming any particular \gls{poc} model: in general, for a given \gls{hbr}, \gls{poc} is expressed as a function of the relative position, and the combined covariance in the B-plane at \gls{tca}
\begin{equation}
    \pc = \pc(\vec{r}_\mathcal{B},\vec{P}_\mathcal{B}).
    \label{eq:generalPoc}
\end{equation}

\subsection{Earth Orbit Dynamics}
To compute the optimal \gls{cam}, starting from the states given by the \gls{cdm}, the maneuverable satellite is back-propagated up to a suitable time to start executing the maneuver. The algorithm is agnostic to the dynamics model utilized, so any representations of the orbital environment can be employed. Since the time frame is relatively short for typical \gls{cam} scenarios, orbital perturbations like high-order gravitational harmonics, atmospheric drag, solar radiation pressure, and third-body attraction do not play a significant role. To showcase the capability of the algorithm, we include the perturbation from the J2 term of the gravitational potential of the planet
\begin{equation}
\begin{split}
    & \dot{\Vec{r}}= \Vec{v} \\
    & \ddot{x} = -\frac{\mu}{r^3} x\left(1 + k_{J_2}\left(1-5z^2/r^2\right) \right) + u_x, \\
    & \ddot{y} = -\frac{\mu}{r^3} y\left(1 + k_{J_2}\left(1-5z^2/r^2\right) \right) + u_y, \\
    & \ddot{z} = -\frac{\mu}{r^3} z\left(1 + k_{J_2}\left(3-5z^2/r^2\right) \right) + u_z, 
\end{split}
    \label{eq:J2dynamics}
\end{equation}
where $\vec{r} = [x \hspace{3pt} y \hspace{3pt} z]\transp$, $\mu\in\mathbb{R}_\mathrm{+}$ is the gravitational constant of the Earth, $k_{J_2}=\frac{3}{2}\left(\frac{R_E}{r}\right)^2 J_2\in\mathbb{R}_\mathrm{+}$, $R_E\in\mathbb{R}_\mathrm{+}$ is the Earth's equatorial radius, $J_2\in\mathbb{R}_\mathrm{+}$ is the Earth's oblateness coefficient, $t\in\mathbb{R}_{[t_0,t_{CA}]}$ is the time domain that goes from an arbitrary starting time $t_0$ to \gls{tca}, and $\vec{u}=[u_x \hspace{3pt} u_y \hspace{3pt} u_z]\transp\in\mathbb{R}^3$ is the acceleration control action. The mass loss due to the maneuver is considered negligible because a \gls{cam} usually involves a small $\dv$ \cite{Hernando-Ayuso2020}. 

\subsection{Cislunar Dynamics}
In the cislunar domain, the motion is described by the \gls{cr3bp} dynamics. The most convenient reference frame to represent the motion of the satellite is the synodic, which rotates with the same angular speed as the orbital motion of the two main bodies, i.e., the Earth and the Moon. The origin of the synodic reference frame is in the barycenter of the Earth-Moon system; the $x$ axis extends from the origin to the Moon's center of mass, the $z$ axis is in the direction of the angular momentum of the system, and the $y$ axis completes the right-handed triad. In this reference frame, the two main bodies are stationary, and the third body typically moves following a chaotic behavior unless it is in a stable orbit.

The physics of the system is nondimensionalized using the characteristic mass $M=6.04564\times10^{15}$ \si{kg}, the characteristic length $D = 384405$ \si{km}, and the characteristic time $T = 375677$ \si{s}. 
The equations of motion can be written in compact form in the synodic frame as\begin{equation}
\begin{split}
& \dot{\vec{r}} = \vec{v} \\
& \ddot{x} = 2 \dot{y} - \Omega_x + u_x, \\
& \ddot{y} = - 2\dot{x} - \Omega_y + u_y, \\
& \ddot{z} = -z - \Omega_z + u_z, 
\end{split}
\label{eq:cr3bpdynamics}
\end{equation}
where $\Omega_x$, $\Omega_y$, and $\Omega_z$ are the partial derivatives of the effective potential \cite{Singh2013}.

\section{CAM OPTIMIZATION PROBLEM}
\label{sec:optProb}
In this section, the \gls{cam} Optimization Problem is set first as an \gls{ocp} and afterward as a \gls{nlp}. 

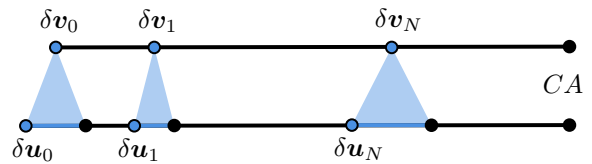
\begin{figure}[b!]
    \centering
    \tikzset{every picture/.style={line width=0.75pt}} %set default line width to 0.75pt        

\begin{tikzpicture}[x=0.75pt,y=0.75pt,yscale=-1,xscale=1]
%uncomment if require: \path (0,465); %set diagram left start at 0, and has height of 465

%Shape: Triangle [id:dp023549451170071234] 
\draw  [draw opacity=0][fill={rgb, 255:red, 74; green, 144; blue, 226 }  ,fill opacity=0.45 ] (190.13,250.41) -- (210.21,289.61) -- (170.06,289.61) -- cycle ;
%Shape: Triangle [id:dp9358994871385596] 
\draw  [draw opacity=0][fill={rgb, 255:red, 74; green, 144; blue, 226 }  ,fill opacity=0.45 ] (70.38,250.41) -- (80.46,289.83) -- (60.31,289.83) -- cycle ;
%Shape: Triangle [id:dp6655658115185338] 
\draw  [draw opacity=0][fill={rgb, 255:red, 74; green, 144; blue, 226 }  ,fill opacity=0.45 ] (20.71,250.4) -- (35.81,289.88) -- (5.6,289.88) -- cycle ;
%Straight Lines [id:da1596466893716737] 
\draw [line width=1.5]    (279.8,289.6) -- (5.6,290) ;
%Straight Lines [id:da32379925391735287] 
\draw [color={rgb, 255:red, 74; green, 144; blue, 226 }  ,draw opacity=1 ][line width=1.5]    (210.21,289.68) -- (170.06,289.61) ;
%Straight Lines [id:da9024194409105077] 
\draw [color={rgb, 255:red, 74; green, 144; blue, 226 }  ,draw opacity=1 ][line width=1.5]    (80.46,289.83) -- (60.31,289.96) ;
%Straight Lines [id:da04871906681492566] 
\draw [color={rgb, 255:red, 74; green, 144; blue, 226 }  ,draw opacity=1 ][line width=1.5]    (35.81,289.88) -- (5.6,290) ;
%Straight Lines [id:da9833112055739309] 
\draw [line width=1.5]    (279.8,250) -- (20.71,250.4) ;
%Shape: Ellipse [id:dp30699175479091567] 
\draw  [fill={rgb, 255:red, 74; green, 144; blue, 226 }  ,fill opacity=1 ] (23.62,250.4) .. controls (23.62,248.79) and (22.31,247.49) .. (20.71,247.49) .. controls (19.1,247.49) and (17.8,248.79) .. (17.8,250.4) .. controls (17.8,252.01) and (19.1,253.31) .. (20.71,253.31) .. controls (22.31,253.31) and (23.62,252.01) .. (23.62,250.4) -- cycle ;
%Shape: Ellipse [id:dp6100330850542428] 
\draw  [fill={rgb, 255:red, 74; green, 144; blue, 226 }  ,fill opacity=1 ] (73.29,250.41) .. controls (73.29,248.8) and (71.99,247.49) .. (70.38,247.49) .. controls (68.78,247.49) and (67.47,248.8) .. (67.47,250.41) .. controls (67.47,252.01) and (68.78,253.32) .. (70.38,253.32) .. controls (71.99,253.32) and (73.29,252.01) .. (73.29,250.41) -- cycle ;
%Shape: Ellipse [id:dp5857136349507337] 
\draw  [fill={rgb, 255:red, 74; green, 144; blue, 226 }  ,fill opacity=1 ] (193.04,250.41) .. controls (193.04,248.8) and (191.74,247.49) .. (190.13,247.49) .. controls (188.53,247.49) and (187.22,248.8) .. (187.22,250.41) .. controls (187.22,252.01) and (188.53,253.32) .. (190.13,253.32) .. controls (191.74,253.32) and (193.04,252.01) .. (193.04,250.41) -- cycle ;
%Shape: Ellipse [id:dp33363392419554116] 
\draw  [fill={rgb, 255:red, 0; green, 0; blue, 0 }  ,fill opacity=1 ] (282.71,250) .. controls (282.71,248.39) and (281.41,247.09) .. (279.8,247.09) .. controls (278.19,247.09) and (276.89,248.39) .. (276.89,250) .. controls (276.89,251.61) and (278.19,252.91) .. (279.8,252.91) .. controls (281.41,252.91) and (282.71,251.61) .. (282.71,250) -- cycle ;
%Shape: Ellipse [id:dp9340240520482856] 
\draw  [fill={rgb, 255:red, 74; green, 144; blue, 226 }  ,fill opacity=1 ] (8.51,290) .. controls (8.51,288.39) and (7.21,287.09) .. (5.6,287.09) .. controls (3.99,287.09) and (2.69,288.39) .. (2.69,290) .. controls (2.69,291.61) and (3.99,292.91) .. (5.6,292.91) .. controls (7.21,292.91) and (8.51,291.61) .. (8.51,290) -- cycle ;
%Shape: Ellipse [id:dp2680120307096274] 
\draw  [fill={rgb, 255:red, 74; green, 144; blue, 226 }  ,fill opacity=1 ] (63.22,289.96) .. controls (63.22,288.35) and (61.91,287.04) .. (60.31,287.04) .. controls (58.7,287.04) and (57.4,288.35) .. (57.4,289.96) .. controls (57.4,291.56) and (58.7,292.87) .. (60.31,292.87) .. controls (61.91,292.87) and (63.22,291.56) .. (63.22,289.96) -- cycle ;
%Shape: Ellipse [id:dp15095608332648558] 
\draw  [fill={rgb, 255:red, 74; green, 144; blue, 226 }  ,fill opacity=1 ] (172.97,289.61) .. controls (172.97,288) and (171.66,286.69) .. (170.06,286.69) .. controls (168.45,286.69) and (167.15,288) .. (167.15,289.61) .. controls (167.15,291.21) and (168.45,292.52) .. (170.06,292.52) .. controls (171.66,292.52) and (172.97,291.21) .. (172.97,289.61) -- cycle ;
%Shape: Ellipse [id:dp7563028030802209] 
\draw  [fill={rgb, 255:red, 0; green, 0; blue, 0 }  ,fill opacity=1 ] (282.71,289.6) .. controls (282.71,287.99) and (281.41,286.69) .. (279.8,286.69) .. controls (278.19,286.69) and (276.89,287.99) .. (276.89,289.6) .. controls (276.89,291.21) and (278.19,292.51) .. (279.8,292.51) .. controls (281.41,292.51) and (282.71,291.21) .. (282.71,289.6) -- cycle ;
%Shape: Circle [id:dp37977584961998323] 
\draw  [fill={rgb, 255:red, 0; green, 0; blue, 0 }  ,fill opacity=1 ] (38.72,289.88) .. controls (38.72,288.27) and (37.42,286.97) .. (35.81,286.97) .. controls (34.2,286.97) and (32.9,288.27) .. (32.9,289.88) .. controls (32.9,291.49) and (34.2,292.79) .. (35.81,292.79) .. controls (37.42,292.79) and (38.72,291.49) .. (38.72,289.88) -- cycle ;
%Shape: Circle [id:dp8204538927987942] 
\draw  [fill={rgb, 255:red, 0; green, 0; blue, 0 }  ,fill opacity=1 ] (83.37,289.83) .. controls (83.37,288.22) and (82.07,286.92) .. (80.46,286.92) .. controls (78.85,286.92) and (77.55,288.22) .. (77.55,289.83) .. controls (77.55,291.44) and (78.85,292.74) .. (80.46,292.74) .. controls (82.07,292.74) and (83.37,291.44) .. (83.37,289.83) -- cycle ;
%Shape: Circle [id:dp6371160273594352] 
\draw  [fill={rgb, 255:red, 0; green, 0; blue, 0 }  ,fill opacity=1 ] (213.12,289.68) .. controls (213.12,288.07) and (211.82,286.77) .. (210.21,286.77) .. controls (208.6,286.77) and (207.3,288.07) .. (207.3,289.68) .. controls (207.3,291.29) and (208.6,292.59) .. (210.21,292.59) .. controls (211.82,292.59) and (213.12,291.29) .. (213.12,289.68) -- cycle ;

% Text Node
\draw (264.82,262.96) node [anchor=north west][inner sep=0.75pt]    {$CA$};
% Text Node
\draw (-2.94,294.79) node [anchor=north west][inner sep=0.75pt]    {$\delta \vec{u}_{0}$};
% Text Node
\draw (50.47,294.74) node [anchor=north west][inner sep=0.75pt]    {$\delta \vec{u}_{1}$};
% Text Node
\draw (159.46,294.49) node [anchor=north west][inner sep=0.75pt]    {$\delta \vec{u}_{N}$};
% Text Node
\draw (180.18,229.49) node [anchor=north west][inner sep=0.75pt]    {$\delta \vec{v}_{N}$};
% Text Node
\draw (60.8,230.05) node [anchor=north west][inner sep=0.75pt]    {$\delta \vec{v}_{1}$};
% Text Node
\draw (11.2,229.69) node [anchor=north west][inner sep=0.75pt]    {$\delta \vec{v}_{0}$};

\end{tikzpicture}
    \caption{Impulsive and low-thrust designs.}
    \label{fig:imp_vs_lt}
\end{figure}

\subsection{Optimal Control Problem}
The fuel-optimal \gls{cam} \gls{ocp} has the objective of finding the minimum control action that can grant a reduction of \gls{poc} below an arbitrary threshold. The thrusting opportunities are fixed at predefined times before \gls{tca}. These are collected into the set $\mathbb{T} = \{t_0,\hspace{3pt} t_1,\hspace{3pt} ...,\hspace{3pt} t_N\}$, $N\in\mathbb{N}$. 
The algorithm proposed in this work is intended to compute maneuvers either using a low-thrust dynamics model or impulsive approximation. In the former case, the optimization variable is the continuous $\vec{u}(t)$. In the latter, $\vec{u}(t) = \vec{0}_3$ at every time, and the control is given by the discrete impulses $\Delta \vec{v}_i$ for $i\in\{0,\hspace{3pt} ...,\hspace{3pt} N\}$. In \cref{fig:imp_vs_lt}, the equivalent impulsive and low-thrust schemes are shown: the blue nodes depict the times where the control variables are placed, while the black ones indicate ballistic propagation points. To include both possibilities, in the following, the control variable is denominated $\vec{\phi}\in\mathbb{R}^3$ and its magnitude $\phi\in\mathbb{R}$.
\begin{subequations}
\begin{align}
    \min_{\vec{u}(t)\lor\Delta\vec{v}_i} \quad & J = \begin{cases}
      \int_{t_0}^{t_N} ||\vec{u}(\tau)||\mathrm{d}\tau \\
      \sum_{i=0}^{i=N} ||\Delta\vec{v}_i||
      \end{cases}
      \label{eq:ocpObj}\\
    \mbox{s.t.} 
    \quad &  \text{\cref{eq:J2dynamics} or \cref{eq:cr3bpdynamics}}  \label{eq:ocpDyn}\\
    \quad & \pc(t) \leq \bar{P}_C \label{eq:ocpPc}\\
    \quad & \vec{r}(t_0) = \vec{r}_0 \label{eq:ocpInit}
    % \quad & ||\vec{u}(t)|| \leq u_{max} \label{eq:ocpUmax},
\end{align}
\label{eq:ocp}
\end{subequations}
\cref{eq:ocpObj} imposes the minimization of the fuel-optimal objective function $J\in\mathbb{R}_\mathrm{+}$: the first case is used for low-thrust dynamics, the second one for the multi-impulsive; \cref{eq:ocpDyn} is the dynamics constraint, which depends on the selected model; \cref{eq:ocpPc} is the \gls{poc} constraint, which imposes that \gls{poc} must be below the required threshold $\pclim$ at any time; \cref{eq:ocpInit} defines the initial position of the satellite, which cannot be altered.

\subsection{Nonlinear Program}
Problem \eqref{eq:ocp} is discretized by using a Runge-Kutta 7-8 integration scheme.
The initial state $\vec{x}_0$ is yielded by back-propagating the ballistic motion of the primary starting from \gls{tca}. The nominal state at subsequent thrusting opportunities can be obtained via successive forward propagations. Therefore, the discretized dynamics equations at any discretization node $t_i$ are generally written
\begin{equation}
\begin{aligned}
    \quad & \vec{x}_i = \vec{f}_i(\vec{x}_{i-1},\vec{u}_{i-1},\vec{p}) \quad & i\in\{1,\hspace{3pt} ...,\hspace{3pt} N\},
\end{aligned}
\end{equation}
where $\vec{f}_i(\cdot):\mathbb{R}^6\times\mathbb{R}^3\times\mathbb{R}^{m_p}\rightarrow\mathbb{R}^6$ is the discretized dynamics function and $\vec{p}\in\mathbb{R}^{m_p}$ is a vector of parameters. The discretization assumes a first-order hold for the input action; therefore, the control between two consecutive firing opportunities is constant.

By means of \gls{da}, we introduce perturbations on the control variable at each thrusting opportunity ($\vec{\phi}_i+\delta\vec{\phi}_i$) in a forward-propagation scheme. Considering a ballistic reference trajectory (no initial control), in the multi-impulsive case, the velocity perturbation is added to the state at every node
\begin{subequations}
\begin{align}
 & \vec{x}_i = \mathcal{T}^n\big|_{\vec{x}_i}(\vec{x}_{i-1}) + [\vec{0}_3 \hspace{3pt} \delta\vec{v}_i]\transp  & i\in\{1,\hspace{3pt} ...,\hspace{3pt} N\} \\
 & \vec{x}_{p,0} = \vec{x}_{p,0} + [\vec{0}_3 \hspace{3pt} \delta\vec{v}_0]\transp,  & 
\label{eq:perturbations1}
\end{align}
\end{subequations}
where $\mathcal{T}^n\big|_{\vec{x}_i}(\cdot):\mathbb{R}^6\rightarrow\mathbb{R}^6$ is the $n^\text{th}$-order Taylor series approximation of $\vec{x}_i$.
In the low-thrust case, the acceleration is defined at every node using \gls{da} perturbations, and it directly contributes to the dynamics\footnote{As illustrated in \cref{fig:imp_vs_lt}, the last node of each low-thrust arc is idle to set where the arc stops. 
This implies that a low-thrust design is more computationally demanding than an impulsive one since at least two nodes are needed to define the low-thrust arc, while only one is necessary for the impulse.}
\begin{equation}
\begin{aligned}
 & \vec{x}_i = \mathcal{T}^n\big|_{\vec{x}_i}(\vec{x}_{i-1},\vec{u}_{i-1})  & i\in\{1,\hspace{3pt} ...,\hspace{3pt} N\}, \\
\end{aligned}
\label{eq:perturbations}
\end{equation}
where $\mathcal{T}^n\big|_{\vec{x}_i}(\cdot):\mathbb{R}^6\times\mathbb{R}^3\rightarrow\mathbb{R}^6$.
A detailed explanation of the use of \gls{da} can be found in \cite{Armellin2010}. As this is a polynomial expansion, its accuracy is inversely proportional to the entity of the perturbation and directly proportional to the polynomial order.

Given the dependence of \gls{poc} on the relative position at \gls{tca}, and the dependence of this from the control history, it is possible to construct a multivariate polynomial representation of \gls{poc} as dependent on the control perturbations stacked vector, defined as $\delta\vec{\Phi} = [\delta\vec{\phi}_0\transp \hspace{3pt} ... \hspace{3pt} \delta\vec{\phi}_N\transp]\transp\in\mathbb{R}^{M}$, where $M = 3(N+1)$
\begin{equation}
   \pc = \mathcal{T}^n\big|_{\pc}(\delta\vec{\Phi}),
    \label{eq:pocPoly}
\end{equation}
where $\mathcal{T}^n\big|_{P_C}(\cdot):\mathbb{R}^M\rightarrow\mathbb{R}$.
We define the vector of the control history as 
\begin{equation}
\vec{\Phi} = \begin{bmatrix} \vec{\phi}_0\transp &  ...  & \vec{\phi}_N\transp\end{bmatrix}\transp\in\mathbb{R}^{M}.
\label{eq:optVect}
\end{equation}
A \gls{nlp} is constructed in the following form
\begin{subequations}
    \begin{align}
        \min_{\Vec{\Phi}} \quad & \sum_{i=0}^{i=N} ||\vec{\phi}_i||\label{eq:objNlp}\\
        \mathrm{s.t.} \quad & \mathcal{T}^n\big|_{\pc}(\vec{\Phi}) = \pclim,
        \label{eq:polConstrNlp}
    \end{align}
    \label{prob:nlp}
\end{subequations}
where \cref{eq:objNlp} is the fuel-optimal objective function and \cref{eq:polConstrNlp} is the scalar polynomial \gls{poc} constraint. The inequality constraint \cref{eq:ocpPc} is turned into an equality sign because it is hypothesized that the \gls{poc} reduction is proportional to the entity of the maneuver.
Problem \eqref{prob:nlp} can be solved using a global optimization tool\footnote{\url{https://yalmip.github.io/tutorial/globaloptimization/}} or MATLAB's fmincon\footnote{\url{https://au.mathworks.com/help/optim/ug/fmincon.html}}.

\section{SEQUENTIAL POLYNOMIAL RECURSIVE METHOD}
\label{sec:recursive}
The recursive method that is proposed in this section is used to solve the energy-optimal counterpart of Problem \eqref{prob:nlp}. Therefore, the \gls{nlp} is converted into a \gls{pp} by changing the objective function, which becomes quadratic
\begin{subequations}
    \begin{align}
        \min_{\Vec{\Phi}} \quad & \Vec{\Phi}\transp\Vec{\Phi}
        \label{eq:obj}\\
        \mathrm{s.t.} \quad & \mathcal{T}^n\big|_{\pc}(\vec{\Phi}) = \pclim.
        \label{eq:polConstr}
    \end{align}
    \label{prob:pp}
\end{subequations}
We employ a recursive approach to find solutions to the \gls{pp} of increasing polynomial orders. The polynomial is first truncated in the $1^\text{st}$ order to obtain an approximate solution; this solution is used as a first guess to solve the problem with a polynomial constraint truncated in the $2^\text{nd}$ order. An iterative process is applied to find the second-order solution; when convergence is reached, the recovered solution is fed to a $3^\text{rd}$-order truncated polynomial, and this process is repeated up to the $n^\text{th}$ order. 

\subsection{Mathematical Background}
\label{sec:background}
The formulation of the recursive method is based on the use of tensor notation for Taylor expansions. 
Given a function $g(\cdot):\mathbb{R}^m\rightarrow\mathbb{R}$, its Taylor expansion can be computed up to an arbitrary order $n$. The polynomial approximating $g(\cdot)$, as computed on the expansion point $\vec{y}^*\in\mathbb{R}^m$ can be written as
\begin{equation}
\mathcal{T}^n\Big|_{g(\vec{y})}(\vec{y}) = g(\vec{y}^*) + \sum_{k=1}^{n} \frac{1}{k!}\frac{\partial^k g}{\partial \vec{y}^k}(\vec{y}^*)(\vec{y}-\vec{y}^*)^k.
\label{eq:TaylorExample}
\end{equation}

\cref{eq:TaylorExample} is equivalent to writing the polynomial using tensors of increasing orders \cite{Jenson2024}. The gradient of $g(\cdot)$ is a $1^\text{st}$-order tensor, i.e., a vector; its Hessian is a $2^\text{nd}$-order tensor, i.e., a matrix, and so on. Therefore, it is possible to define the $k$\textsuperscript{th}-order tensor using $k$ different indices: the vector $\vec{a}\in\mathbb{R}^m$ is indexed as $a_i$ ($i\in\{1,\hspace{3pt} ...,\hspace{2pt} m\}$), the matrix $\vec{A}\in\mathbb{R}^{m^2}$ is indexed as $A_{i_1 i_2}$ ($i_1,i_2\in\{1,\hspace{3pt} ...,\hspace{2pt} m\}$), the third order tensor $\vec{\mathcal{A}}\in\mathbb{R}^{m^3}$ is indexed $\mathcal{A}_{i_1 i_2 i_3}$ ($i_1,i_2,i_3\in\{1,\hspace{3pt} ...,\hspace{2pt} m\}$), and so on. With this tool, we can now express the $k$\textsuperscript{th}-order derivatives of $g(\cdot)$ using $k$\textsuperscript{th}-order tensors. Moreover, since the derivative operation is commutative, these tensors have the useful property of super-symmetry, i.e., they are invariant to permutations of the indices. In the following, the order of the tensor will be indicated as a bracketed superscript, $\vec{\mathcal{A}}^{(k)}$. Given a tensor $\vec{\mathcal{A}}^{(k)}$ and a vector $\vec{v}\in\mathbb{R}^m$, their multi-linear form  $\vec{\mathcal{A}}^{(k)} \vec{v}^k \in\mathbb{R}$ is defined as:
\begin{equation}
   \vec{\mathcal{A}}^{(k)} \vec{v}^k = \sum_{i_1,i_2,...,i_k = 1}^m \mathcal{A}_{i_1i_2...i_k} v_{i_1}v_{i_2}...v_{i_k},
\label{eq::multilinear}
\end{equation}
where the notation $\vec{v}^k$ refers to the repetition of $\vec{v}$  for $k$ times. The summation is only performed over indices that are unique to the right-hand side of the equation.
With this knowledge, we can now re-write \cref{eq:TaylorExample} using multi-linear forms
\begin{equation}
\mathcal{T}^n\Big|_{g(\vec{y})}(\vec{y}) = g(\vec{y}^*) + \sum_{k=1}^{n} \vec{\mathcal{G}}^{(k)} (\vec{y}-\vec{y}^*)^k,
\label{eq:TaylorExample1}
\end{equation}
where $\vec{\mathcal{G}}^{(k)}$ is the $k$\textsuperscript{th}-order tensor that collects the $k$\textsuperscript{th}-order derivatives with respect to $\vec{y}$.

If the multiplication by the first tensor mode is eliminated, we obtain a row-vector output\footnote{In the literature it is often assumed to be a column-vector, but a row-vector is preferred for the purpose of this work since it will be used to define a pseudo-gradient.}, indicated as $\vec{\mathcal{A}}^{(k)} \vec{v}^{k-1}\in\mathbb{R}^m$. Its definition is as follows
\begin{equation}
   (\vec{\mathcal{A}}^{(k)} \vec{v}^{k-1})_j = \sum_{i_2,...,i_k = 1}
    ^m\mathcal{A}_{j i_2 i_3 ...i_k} v_{i_2} v_{i_3}...v_{i_k},
\label{eq:noFirstForm}
\end{equation}
which is valid for each element of the row $j\in\{1,\hspace{3pt} ...,\hspace{2pt} m\}$.

\subsection{Method}

Let us re-write the polynomial $\mathcal{T}^n\big|_{\pc}(\vec{\Phi})$ from \cref{eq:polConstr} using the multi-linear forms introduced in \cref{sec:background}
\begin{subequations}
\begin{equation}
    \mathcal{T}^n\big|_{\pc}(\vec{\Phi}) = P_C^0 + \sum_{k=1}^n \vec{\mathcal{F}}^{(k)} \vec{\Phi}^k,
    \label{eq:polyExplicit}
\end{equation}
\begin{equation}
    \vec{\mathcal{F}}^{(k)} = \frac{1}{k!}\frac{\partial^k P_C}{\partial \vec{\Phi}^k}
    \label{eq:tensorsDefinition}
\end{equation}
\begin{equation}
    P_C^0 = \mathcal{T}^0\Big|_{P_C}(\vec{\Phi}) 
\end{equation}
\end{subequations}
where $P_C^0\in\mathbb{R}$ is \gls{poc} of the unperturbed trajectory, i.e., the $0^\text{th}$-order of the Taylor polynomial; $\vec{\mathcal{F}}^{(k)}\in\mathbb{R}^{M^k}$ - for $k\in\{1,\hspace{3pt} ...,\hspace{2pt} n\}$ - are symmetric $k$\textsuperscript{th}-order tensors that represent the contributions of all the $k$\textsuperscript{th}-order partial derivatives of \gls{poc} with respect to $\vec{\Phi}$. The number of dimensions of the tensor $\vec{\mathcal{F}}^{(k)}$ is equal to the order of the associated polynomial term, e.g., $\vec{\mathcal{F}}^{(2)}\in\mathbb{R}^{M^2}$, $\vec{\mathcal{F}}^{(3)}\in\mathbb{R}^{M^3}$, and so on.

The first-order truncation of \cref{eq:polyExplicit} reads
\begin{equation}
     \mathcal{T}^1\big|_{\pc}(\vec{\Phi}) = P_C^0 + (\vec{\nabla}P_C)\vec{\Phi},  
     \label{eq:first-order1}
\end{equation}
where $\vec{\nabla}P_C=\vec{\mathcal{F}}^{(1)}\in\mathbb{R}^M$ is the gradient of \gls{poc} with respect to the stacked control (row vector).
The probability gap that needs to be filled by the maneuver is called $\rho\in\mathbb{R}$
\begin{equation}
    \rho := \bar{P}_C - P_C^0. 
    \label{eq:definitionRho}
\end{equation}
Then, the first-order truncation of the equality constraint \cref{eq:polConstr} becomes
\begin{equation}
    \rho = (\vec{\nabla}P_C)\vec{\Phi}.
\end{equation}
This equation suggests that, in a first-order approximation, a thrust in the direction of the gradient grants the highest possible change in \gls{poc}. Therefore, we compute the greedy solution of the first-order constraint as
\begin{equation}
    [\vec{\Phi}]_1 = \frac{\rho}{\nabla P_C}\hat{\nabla} P_C,
     \label{eq:first-order}
\end{equation}
where $\nabla  P_C$ is the norm of the gradient, $\hat{\nabla} P_C$ is its direction, and the index $1$ indicates that this is a solution to the first-order constraint.

Now, let us assume that a greedy solution has been found for the $(j-1)$\textsuperscript{th}-order polynomial program; the generalized process for finding the $j$\textsuperscript{th}-order solution is as follows.
The $j$\textsuperscript{th}-order truncation of \cref{eq:polConstr} is written using the multi-linear form notation from \cref{eq:polyExplicit} and the definition in \cref{eq:definitionRho}

\begin{equation}
     \rho = \sum_{k=1}^j \vec{\mathcal{F}}^{(k)} \vec{\Phi}^k.
    \label{eq:j-order}
\end{equation}
In this case, finding a greedy solution is not as simple as in the first-order case because we have a higher-order dependence from $\vec{\Phi}$. Therefore, successive linearizations of \cref{eq:j-order} are carried on until convergence. The iteration number is indicated with $b\in\mathbb{N}$. In the first iteration, the linearization point is the output of the previous order $\tilde{\vec{\Phi}} = [\vec{\Phi}]_{j-1}^{b_{end}}$. To express the linearized polynomial, we make use of the multi-linear mode with the exclusion of the first tensor mode that was defined in \cref{eq:noFirstForm}
\begin{equation}
     \rho = \sum_{k=1}^j \left(\vec{\mathcal{F}}^{(k)} \tilde{\vec{\Phi}}^{k-1}\right)\vec{\Phi}.
     \label{eq:j-orderLinear}
\end{equation}
\cref{eq:j-orderLinear} is a linear function because $\tilde{\vec{\Phi}}$ is a known value. Therefore, it is possible to define a pseudo-gradient which includes the contribution of the high-order terms  
\begin{equation}
    \vec{g}_j = \vec{\nabla}P_C + \sum_{k=2}^j \vec{\mathcal{F}}^{(k)} \tilde{\vec{\Phi}}^{k-1},
    \label{eq:pseudo-gradient}
\end{equation}
Now, one can compute the greedy solution of the 
$j$\textsuperscript{th}-order constraint
\begin{equation}
    [\vec{\Phi}]_j^b = \frac{\rho}{g_j}\hat{g}_j.
    \label{eq:greedyj}
\end{equation}
With this solution, the new linearization point is selected
\begin{equation}
    \tilde{\vec{\Phi}} = [\vec{\Phi}]_j^b.
    \label{eq:newY}
\end{equation}
Convergence is checked by evaluating the norm of the difference between the solution of iteration $b$ and the previous one, $b-1$ 
\begin{equation}
    || [\vec{\Phi}]_j^b - [\vec{\Phi}]_j^{b-1} ||_2 \leq e_{tol},
    \label{eq:convj}
\end{equation}
where $e_{tol}\in\mathbb{R}$ is an arbitrarily small number.
If \cref{eq:convj} is not satisfied, the procedure is repeated from \cref{eq:j-orderLinear} to \cref{eq:newY}. When convergence within the iterations is reached, the truncation order is increased up to the original order of the constraint $n$.

\begin{figure}[b!]
    \centering
    \input{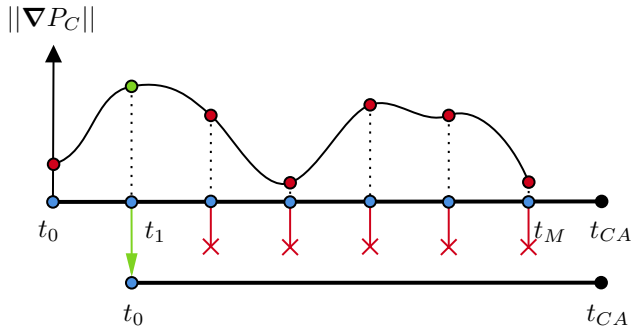}
    \caption{Maneuvering times ranking and selection routine.}
    \label{fig:gridRank}
\end{figure}

The flexibility of the method can be exploited to design operation-oriented routines. A nodes-filtering approach is based on the fact that we can rank the entity of the effect a maneuver performed at different times has on \gls{poc} using the gradient norm. A dense grid of available thrusting opportunities is built between $t_0$ and $t_{CA}$, which is represented by the set $\{t_0, \hspace{3pt} ..., \hspace{3pt} t_M\}$. First-order maps of \gls{poc} are computed with respect to an impulse (or thrust arc) at each point on the grid. If $N\in\mathbb{N}$ is the number of thrust opportunities we want to use, only the first $N$ nodes where the gradient is maximum are kept, while the others are eliminated from the discretization. For example, in \cref{fig:gridRank}, out of the $M$ available maneuvering times, only the second one is kept, and the optimization is run using a single impulse. From this point onward, the algorithm is run as presented in the previous paragraphs using the filtered nodes. 

If the maximum thrust is limited, the nodes-filtering routine can be further taken advantage of. The nodes are ranked as explained in the previous paragraph, and, initially, only the first node is kept; if the magnitude of the optimized impulse exceeds the limit imposed by the thrusters, a new reference trajectory is computed subject to a thrust in the optimized direction. A new optimization is run with the inclusion of the second-ranked node, considering a reduced \gls{poc} gap ($\rho$) due to the influence of the first insufficient impulse. The process is repeated if necessary until the \gls{poc} gap is reduced to $0$.

The method can also be used to design optimal maneuvers when the control direction is fixed. This case is common for operators who prefer not to change the attitude of the spacecraft and typically fire in the tangential direction when performing a \gls{cam}; in this case, the optimization vector $\vec{\Phi}\in\mathbb{R}^{N+1}$ includes the magnitude of the impulse (or the continuous thrust) for each thrusting opportunity. 

\section{RESULTS}
\label{sec:results}

\begin{figure}[b!]
    \centering
    \includegraphics[width=\columnwidth]{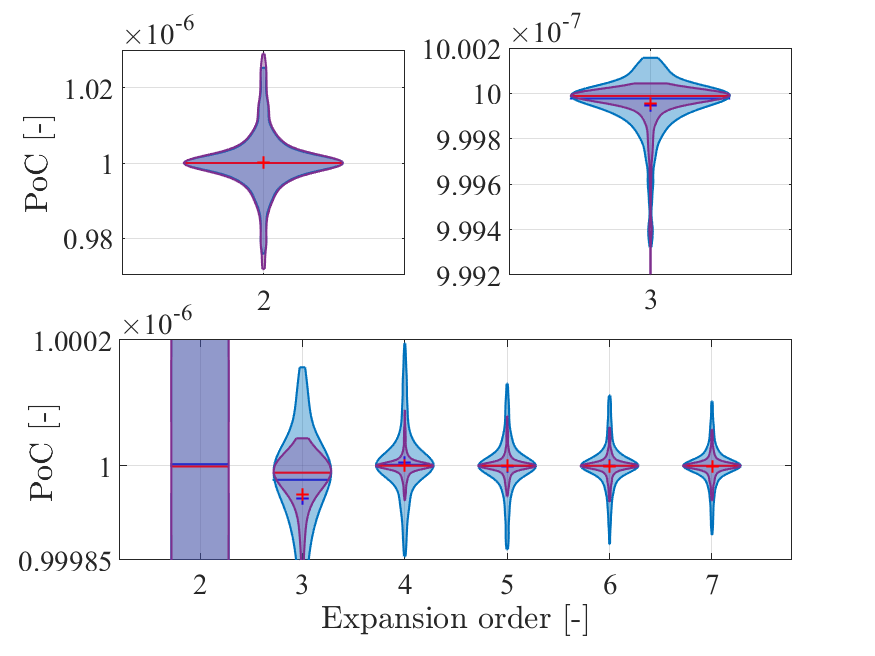}
    \vspace{-.6cm}
    \caption{Distribution of the targeted \gls{poc} in the single impulse campaign for the recursive method (blu violins) and fmincon (purple violins).}
    \label{fig:distributionPoc}
\end{figure}

\begin{figure}[b!]
    \centering
    \includegraphics[width = \columnwidth]{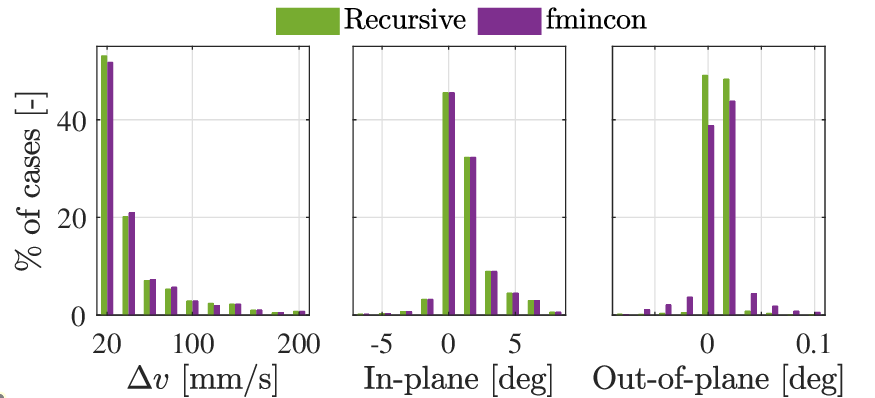}
    \vspace{-.4cm}
    \caption{Distribution of the $\dv$ magnitude and direction for the single-impulse simulation campaign with a $5^\text{th}$-order expansion.}
    \label{fig:dvHistograms}
\end{figure}

\begin{figure}[b!]
    \centering
    \includegraphics[width = \columnwidth]{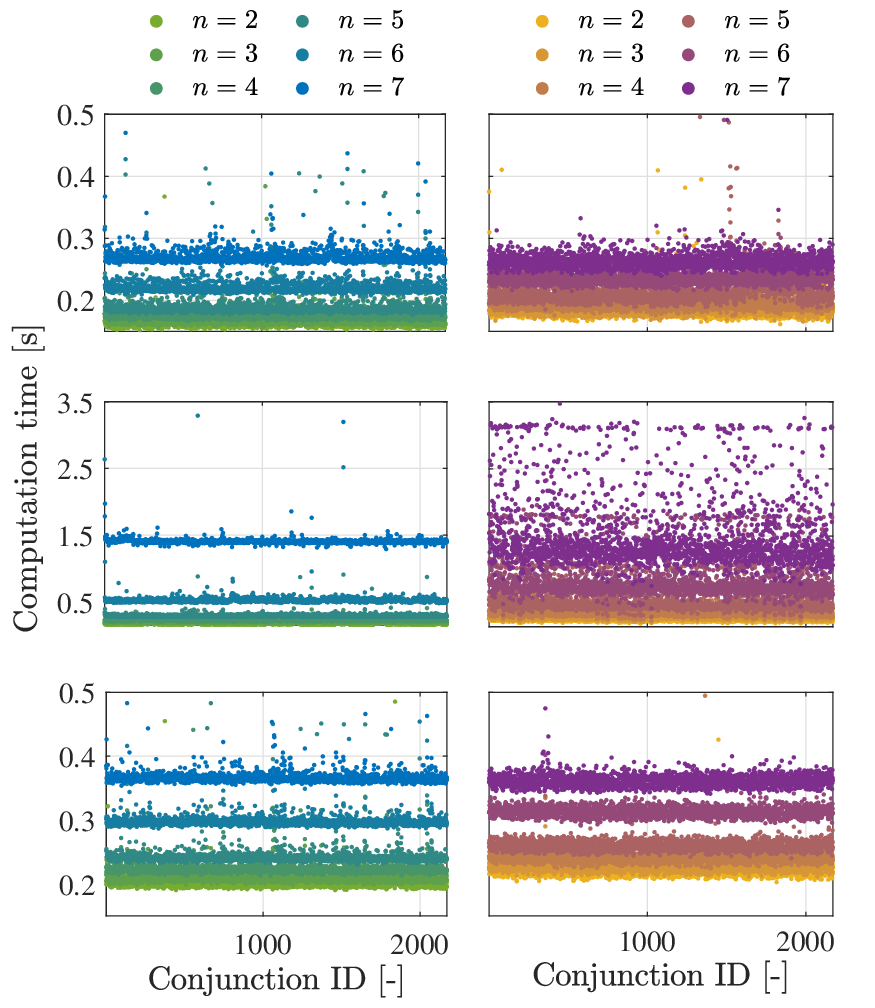}
    \vspace{-.4cm}
    \caption{Computation time for the recursive method (left) and fmincon (right). First row: single-impulse campaign. Second row: bi-impulsive campaign. Third-row: J2 campaign.}
    \label{fig:compTime}
\end{figure}

The results shown in this section include test cases for \gls{leo} and cislunar dynamics. The \gls{leo} ones are taken from the ESA Collision Avoidance Challenge, presented in reference \cite{Armellin2021}, which comprises 2170 scenarios. In all the following analyses, $e_{tol} = 10^{-10}$, the target \gls{poc} is $10^{-6}$, and the dynamics are Keplerian unless otherwise specified. The control components are in \gls{rtn} for Earth orbit scenarios and synodic for the Cislunar one. All the simulations are run on an 11th Gen Intel(R) Core(TM) i7-11700 @2.5GHz, 2496 Mhz, 8 Cores, 16 Logical Processors, using MATLAB R2021b and C++.

\subsection{Extensive LEO campaign}
\label{sec:esacampaign}
First, the recursive method is compared in terms of accuracy, i.e., final targeted \gls{poc}, and computation time, with fmincon's interior-point, used to solve Problem~\eqref{prob:nlp}. fmincon always uses the ballistic trajectory as a first guess. The recovered solution is fed to a forward propagation scheme to validate the results. fmincon and the recursive method solve the \gls{nlp} in Problem \eqref{prob:nlp} and the \gls{pp} in Problem \eqref{prob:pp}, respectively, so fmincon is expected to provide better results thanks to the fuel-optimal formulation. However, it’s worth noticing that the fuel- and energy-optimal formulations coincide in cases of a single impulse or burn. Thus, in the latter case, comparing the two methods provides information about the capability of the recursive method to compute an optimal solution. On the other hand, we can also observe the optimality loss due to the energy-optimal formulation for cases with multiple thrust opportunities. In \cref{fig:distributionPoc}, the distribution of validated \gls{poc} after a maneuver performed $2.5$ orbits before the encounter is shown for different expansion orders. In the violin plots, the recursive method is represented by the light blue areas and fmincon by the purple ones. The mean and median values of the distributions are shown as a cross and a continuous line, respectively: the blue features are relative to the recursive method, the red ones to fmincon. The violin plots show the $98 \%$  of the data. Evidently, the $2^\text{nd}$- and $3^\text{rd}$-order expansions present a larger variance than the following ones, and fmincon can target a \gls{poc} with a higher accuracy. The recursive method becomes very accurate from the $5^\text{th}$-order, showing small improvements for the successive ones. The $\dv$ distribution of the two methods is very similar, as shown in \cref{fig:dvHistograms}: the preferred direction of most test cases is almost tangential, i.e., with null in-plane and out-of-plane angles, and fmincon's solution typically has a slightly higher out-of-plane component. The computation time, shown in \cref{fig:compTime}, is typically very similar for the two methods and proportional to the expansion order. The computation of the Taylor expansion \cref{eq:pocPoly} takes up most of the total run time ($99\%$ at $2^\text{nd}$-order, $85\%$ at $7^\text{th}$-order). On average, among the $2170$ test cases, the number of iterations used by the recursive method with $2^\text{nd}$ to $7^\text{th}$-orders are $60$, $31$, $13$, $6$, $5$, and $5$, respectively. 

The full dataset of conjunctions is also addressed with a bi-impulsive strategy: the first impulse is at $2.5$ orbits before \gls{tca}, the second one at $0.5$. The final \gls{poc} distribution is very similar to the ones reported in \cref{fig:distributionPoc}. The second row of \cref{fig:compTime} shows that, while for lower expansion orders, the recursive method is only slightly faster, above the $5^\text{th}$ order, the computations are significantly faster. In particular, the $5^\text{th}$-order expansion is confirmed as a good choice since its run time is almost always below $0.3$~\si{s}. Since the number of control variables is doubled, in this case the optimization takes a larger portion of the total run-time: for the $2^\text{nd}$-order the recursive method uses $2\%$ of the run-time, while fmincon uses $35\%$; at the $7^\text{th}$-order the percentages go up to $53\%$ for the recursive and $60\%$ for fmincon. In this case, on average, fewer iterations are needed to converge, namely, from $2^\text{nd}$ to $7^\text{th}$-orders, $51$, $22$, $6$, $3$, $2$, and $1$.
The $\dv$ of the bi-impulsive strategy has a slightly different distribution than the single-impulse counterpart: the recursive method finds slightly higher $\dv$ than fmincon, as it solves an energy-optimal problem rather than a fuel-optimal one. As a result, fmincon typically finds a solution that is close to a single impulse, while the recursive method always uses both opportunities to fire. If the distributions of the sum of the two $\dv$ are fitted using a Rayleigh curve, the scale parameter of the recursive method is $77.8$~\si{mm/s}, while fmincon's is $74.2$~\si{mm/s}.

A last campaign is performed using a more complex dynamics model, with the inclusion of the J2 perturbation. In this case, the final solutions have \gls{poc} and $\dv$ distributions that are very similar to the ones in \cref{fig:distributionPoc} and \cref{fig:dvHistograms}. As expected, over the course of a few orbits, the inclusion of high-order harmonics does not influence the system's dynamics appreciably. The run-times, shown in the third row of \cref{fig:compTime}, are generally comparable between the two methods. The increase in the computation time for increasing orders is mainly dictated by the complexity of the \gls{da} propagation. Once the coefficients are obtained, the two methods perform similarly.

\subsection{Single impulse at different times}

\begin{figure}[tb!]
    \centering
    \includegraphics[width = \columnwidth]{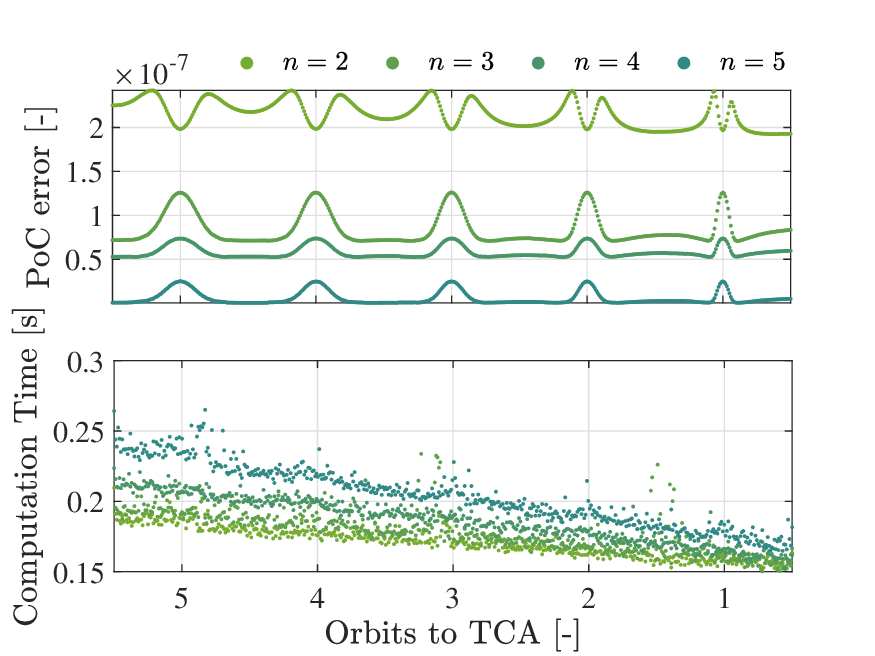}
    \vspace{-.5cm}
    \caption{Targeted \gls{poc} and required computation time of the recursive method for different maneuvering times before \gls{tca}.}
    \label{fig:pocSimTimeDiffManTimes}
\end{figure}

\begin{figure}[tb!]
    \centering
    \includegraphics[width=\columnwidth]{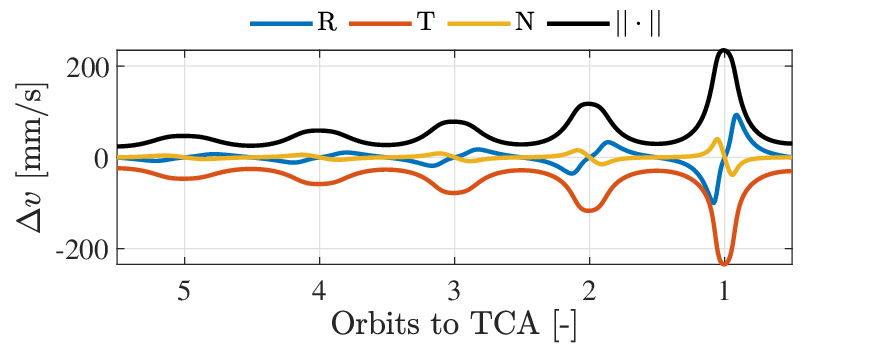}
    \vspace{-.6cm}
    \caption{Evolution of the required $\dv$ for different maneuvering times before \gls{tca} using a $5^\text{th}$-order expansion.}
    \label{fig:dvLeo}
\end{figure}

\begin{figure}[tb!]
    \centering
    \includegraphics[width=\columnwidth]{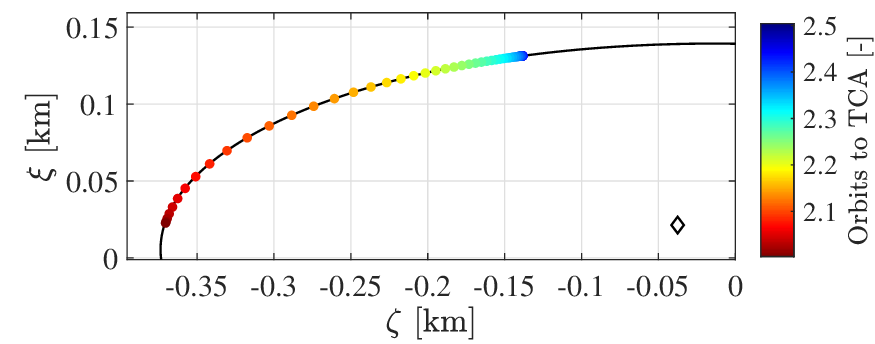}
    \vspace{-.5cm}
    \caption{B-plane configuration in conjunction \#1 for different maneuvering times using a $5^\text{th}$-order expansion. The diamond is the ballistic relative position at \gls{tca}, and the scattered points indicate the position after the maneuver.}
    \label{fig:bplane}
\end{figure}

A detailed analysis of conjunction \#1 is performed to determine the optimal firing time between $0.5$ and $5.5$ orbits before the conjunction. This conjunction represents one of the outliers that were excluded from \cref{fig:distributionPoc}. Indeed, in the first plot of \cref{fig:pocSimTimeDiffManTimes}, one may notice that the targeted \gls{poc} error is quite high for low expansion orders. For the $5^\text{th}$-order, \gls{poc} is very precise at the half-orbits, but it becomes less accurate at the full orbits. In the full orbits, moreover, the required $\dv$ is much higher than in the half-orbits, as shown in \cref{fig:dvLeo}. \cref{fig:bplane} shows that a maneuver close to the half-orbit allows the spacecraft to move along the semi-minor axis of the iso-probability ellipse, where the gradient of \gls{poc} in the B-plane is higher \cite{Pavanello2023Long}. At the full orbits, however, the maneuver causes a shift along the ellipse's semi-major axis, where \gls{poc} changes less rapidly.
For this reason, in agreement with the literature \cite{Bombardelli2015,Armellin2021,DeVittori2023Leo,Pavanello2024Multiple}, we find that firing at half-orbits is best for fuel efficiency. From the second plot in \cref{fig:pocSimTimeDiffManTimes}, it is confirmed that the $5^\text{th}$-order expansion offers the optimal balance between accuracy and computation speed.

\subsection{Low-thrust maneuver}

\begin{figure}[tb!]
    \centering
    \includegraphics[width=\columnwidth]{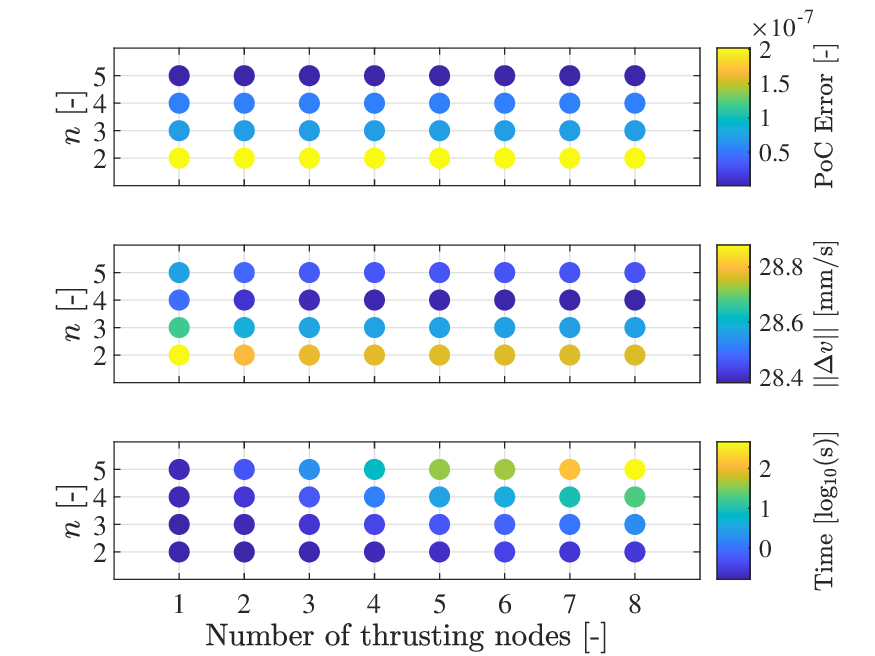}
    \vspace{-.5cm}
    \caption{Accuracy, required $\dv$, and computation time for a single low-thrust window with different discretization points and expansion orders.}
    \label{fig:OrdVsImp}
\end{figure}

\begin{figure}[b!]
    \centering
    \includegraphics[width=\columnwidth]{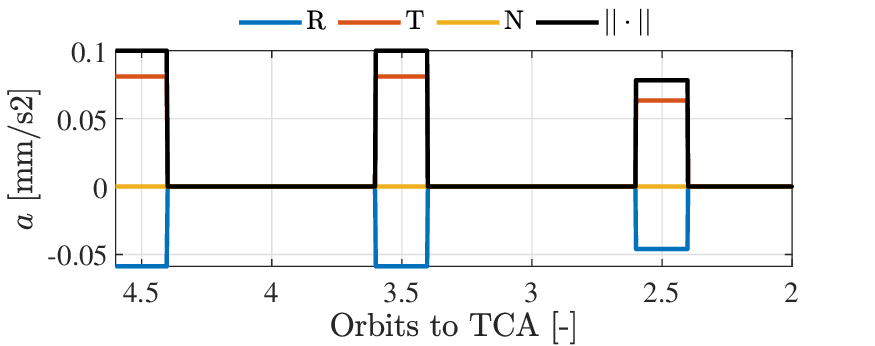}
    \vspace{-.5cm}
    \caption{Acceleration profile using the nodes-filtering routine in conjunction \#1219.}
    \label{fig:accThreeWind}
\end{figure}

The analysis performed on the whole set of conjunctions using the low-thrust model with a single firing opportunity is similar to the impulsive case and is not reported for brevity. Instead, we report an analysis encompassing a different discretization of the thrusting window using different expansion orders. The considered test case is still conjunction \#1. The thrusting window is centered around $2.5$ orbits before \gls{tca}, spans $6$~\si{minutes}, and is discretized in a number of nodes that varies from 1 to 8. The denser the discretization grid is, the more accurate the maneuver can be in terms of thrusting direction. The first plot in \cref{fig:OrdVsImp} shows that the accuracy of the maneuver is independent of the number of discretization nodes. Moreover, the sensitivity of the total required $\dv$ to the number of nodes is very low, too. Given that the computation time, shown in the last plot, grows proportionally to the two variables, a low number of discretization nodes is preferable in this case.

Let us now consider conjunction \#1219, for which the required $\dv$ from the campaign in \cref{sec:esacampaign} is high ($342.5$~\si{mm/s}). The maximum acceleration of a $500$~\si{kg} satellite mounting a $50$~\si{mN} thruster is $0.1$~\si{mm/s^2}. This means that to achieve $342.5$~\si{mm/s} one should thrust for almost $1$~\si{hour}, which is more than half the orbital period. So, the maneuver must be split into multiple thrust arcs over multiple orbits. Given a series of available thrust arcs, centered at every half-orbit from 5 orbits before \gls{tca} and lasting for $20$~\si{minutes}, we can apply the nodes-filtering routine to optimize the maneuver. The resulting thrust profile, shown in \cref{fig:accThreeWind}, grants a targeting of \gls{poc} with an absolute precision of $10^{-10}$: only the first three thrust windows are used. This is not a bang-bang profile, but it can be easily regularized by slightly reducing the thrust time in the last window.
\subsection{Fixed direction}

\begin{figure}[tb!]
    \centering
    \includegraphics[width=\columnwidth]{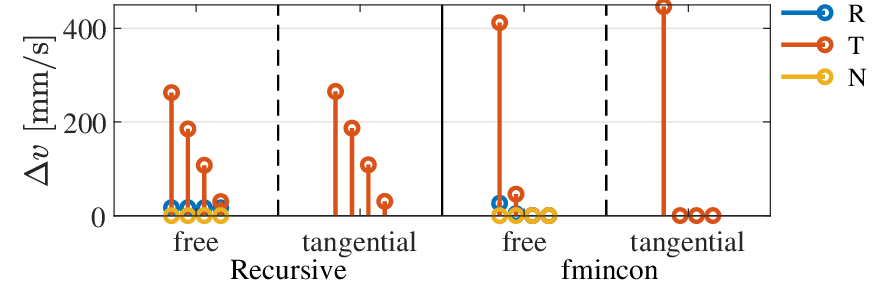}
    \vspace{-.5cm}
    \caption{Tetra-impulsive optimization comparison between the recursive method and fmincon with and without fixed tangential direction.}
    \label{fig:dvTangVsfree}
\end{figure}

Space operators might prefer to thrust in a fixed direction to simplify attitude control. Knowing that when sufficient warning time is given, the tangential direction is the most efficient for \gls{cam}, working with a fixed tangential thrust might simplify the maneuver implementation with little sacrifice in optimality. Additionally, fixing the direction reduces the problem variables to 1/3 of an unconstrained problem, thus improving numerical efficiency.
To exemplify this, let us consider conjunction \#1,466 from the dataset, which is one of the most demanding. Keeping the $5^\text{th}$-order expansion, we compare performances of fmincon and the recursive method using free-direction optimization (optimizing the three thrust components) or tangential impulses. Namely, impulses are available $3.5$, $2.5$, $1.5$ and $0.5$ orbital periods before \gls{tca}. As can be seen from \cref{fig:dvTangVsfree}, when multiple impulses are available, the recursive method finds a solution that uses every available opportunity to fire because the gradient of \gls{poc} is non-null with respect to each thrusting window. On the contrary, fmincon, solving the fuel-optimal problem, uses as few impulses as possible: in the tangential firing case, only one impulse is used. In every case, the validated final \gls{poc} value is achieved with an absolute precision of $10^{-10}$. While obtaining the same accuracy and finding a similar solution, the fixed-direction optimization is much faster than the free-direction one, keeping the execution time below $0.4$~\si{s} both using fmincon ($0.38$~\si{s} vs $3.41$~\si{s}) and the recursive method ($0.26$~\si{s} vs $0.86$~\si{s}). This significant increase in efficiency comes at a very low loss in terms of optimality: the total $\dv$ is only increased by $1.2$~\si{mm/s} using the recursive method and $12.6$~\si{mm/s} with fmincon.

\subsection{Cislunar dynamics}
\begin{figure}[tb!]
    \centering
    \includegraphics[width=\columnwidth]{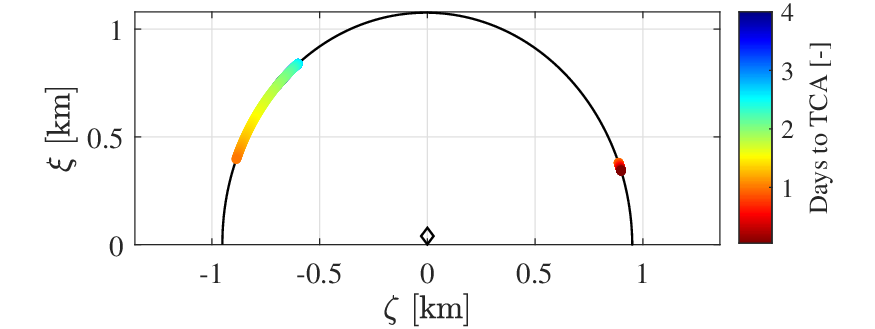}
    \vspace{-.5cm}
    \caption{B-plane configuration in the cislunar test-case for different maneuvering times using a $2^\text{nd}$-order expansion. The blue kernels are hidden behind the light blue ones.}
    \label{fig:bplaneCislunar}
\end{figure}

\begin{figure}[tb!]
    \centering
    \includegraphics[width=\columnwidth]{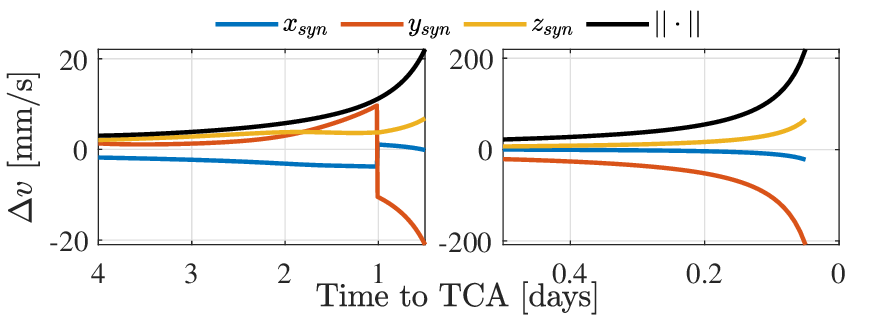}
    \vspace{-.6cm}
    \caption{Evolution of the required $\dv$ for different maneuvering times before \gls{tca} using a $2^\text{nd}$-order expansion in the cislunar test case.}
    \label{fig:dvCislunar}
\end{figure}

The single-impulse optimization is applied to a L1-Near
Rectilinear Halo Orbit from reference \cite{DeMaria2023} to showcase the applicability to the cislunar domain. The $2^\text{nd}$-order expansion is accurate enough in this application, reaching very similar values of final \gls{poc} compared to higher orders. In \cref{fig:bplaneCislunar}, the positions reached on the B-plane for maneuvering times spanning from $4$ days to $1$ hour before \gls{tca} are shown. In \cref{fig:dvCislunar}, the required impulsive $\dv$ components are presented. The results are very similar to the ones obtained with the analytical methods from \cite{DeMaria2023}. In particular, both algorithms find two different optimal wells that are targeted alternatively depending on the alert time: in the first three days, the region in the second quadrant of the B-plane is targeted, while for shorter alert times, the region in the first quadrant is preferred. The computational time is comparable to the \gls{leo} test cases: out of the $500$ simulations performed to analyze the test case, the mean execution time was $0.153$~\si{s}, with a standard deviation of $0.005$~\si{s}.

\section{CONCLUSIONS}
\label{sec:conclusions}
A recursive method for the design of a \acrfull{cam} was presented based on the polynomial expansion of the \acrfull{poc}. The method is generally applicable to scenarios involving a single short-term encounter with any dynamics, and it can deal with an arbitrary number of different maneuvering windows. Moreover, it can deal both with low-thrust and with impulsive dynamics, optimizing the control acceleration of the thrust arc in the first case and the impulses in the second. 

The computational results showcase various combinations of \gls{cam} dynamics, thrust models, maneuvering times, and expansion orders. The recursive method's solution to the energy-optimal \acrlong{pp} is compared with fmincon's, which is a solution to the fuel optimal \acrlong{nlp}. Very similar results are obtained by the two methods in most of the cases: a target \gls{poc} is almost always achieved with very high accuracy, and the computational time is usually kept below $1$~\si{s}. The simplicity and immediateness of the method, combined with the quasi-optimal results that it can achieve, make it potentially suitable for onboard implementation.

\begin{footnotesize}
\bibliography{references}
\end{footnotesize}
\begin{IEEEbiography}[{\includegraphics[width=1in,height=1.25in,clip,keepaspectratio]{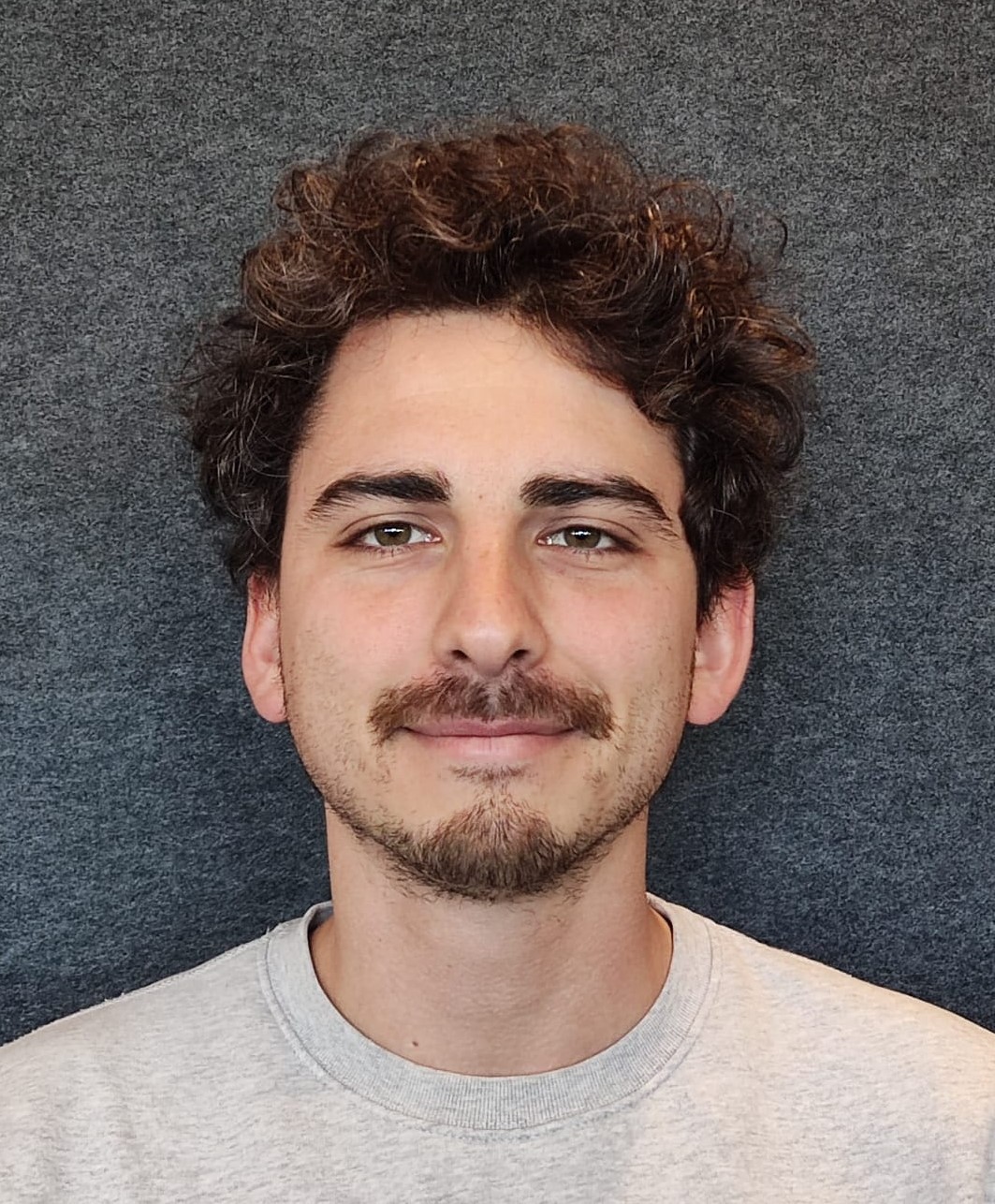}}]
{Zeno Pavanello}{\space}
obtained his B.Sc. degree in Aerospace Engineering from the University of Padua in 2017. He obtained his M.Sc. double degree in Aerospace Engineering in 2020 from the University of Padua and the Instituto Superior Técnico of Lisbon. Zeno is currently a third-year Ph.D. student at Te P\=unaha \=Atea – Space Institute at the University of Auckland. He works on collision avoidance maneuvers optimization.
\end{IEEEbiography}
\vspace{-0.4cm}
\begin{IEEEbiography}[{\includegraphics[width=1in,height=1.25in,clip,keepaspectratio]{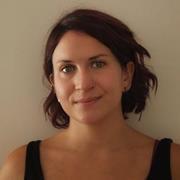}}]{Laura Pirovano} received her M.Sc. degree in aerospace engineering (space exploration track) from TU Delft, the
Netherlands, in 2015. She completed her Ph.D. at the University of Surrey, UK, in 2020 on methods for cataloging space debris with optical observations. She is now a research fellow at Te P\=unaha \=Atea – Space Institute at the University of Auckland. Her research interests include space situational awareness and uncertainty propagation.
\end{IEEEbiography}
\vspace{-0.4cm}
\begin{IEEEbiography}[{\includegraphics[width=1in,height=1.25in,clip,keepaspectratio]{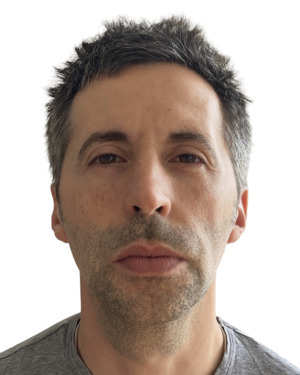}}]{Roberto Armellin} received his M.Sc. and Ph.D. degrees in aerospace engineering from Politecnico di Milano, Italy, in 2003 and 2007, respectively. Since November 2020, he has been a professor at Te P\=unaha \=Atea – Space Institute, the University of Auckland. His current research interests include space trajectory optimization, spacecraft navigation and guidance, and space situational awareness. 
\end{IEEEbiography}
\end{document}